\documentclass[12pt]{amsproc}
\usepackage{amsmath,amssymb,amscd}

\begin{document}

\title[Free Probability for Pairs of Faces IV]{Free Probability for Pairs of Faces IV: \\
Bi-Free Extremes in the Plane}
\author[Dan-Virgil Voiculescu]{Dan-Virgil Voiculescu \\ Department of Mathematics \\ University of California at Berkeley \\ Berkeley, CA\ \ 94720-3840 \\ E-mail: dvv@math.berkeley.edu}
\address{}
\thanks{Research supported in part by NSF Grant DMS-1301727.}
\keywords{Bi-free independence, Bi-free extremal convolution, Two-faced pairs.}
\subjclass[2000]{Primary: 46L54; Secondary: 60G70, 46L53}
\date{}

\begin{abstract}
We compute the bi-free max-convolution which is the operation on bi-variate distribution functions corresponding to the max-operation with respect to the spectral order on bi-free bi-partite two-faced pairs of hermitian non-commutative random variables. With the corresponding definitions of bi-free max-stable and max-infinitely-divisible laws their determination becomes in this way a classical analysis question.
\end{abstract}

\maketitle

\setcounter{section}{-1}
\section{Introduction}
\label{sec0}

The definition and classification of free max-stable laws in \cite{2} had been an unexpected addition to the list of free probability analogues to classical probability items. Here we take the first step in a similar direction in bi-free probability \cite{9}. We show that there is a simple formula for computing the bi-free extremal convolution of probability measures in the plane. This corresponds to computing the distribution of $(a \vee c,b \vee d)$, where $(a,b)$ and $(c,d)$ are two bi-free two-faced pairs of commuting hermitian operators. Like the free extremal convolution on ${\mathbb R}$ defined in \cite{2}, the bi-free extremal convolution in the plane reduces the question of bi-free max-stable laws in the plane to an analysis problem in the classical context which we won't pursue here.

The derivation of the formula for the extremal bi-free convolution relies on the partial bi-free $R$-transform we found in \cite{10}. Obviously our result here complements the recent work on operations on bi-free bi-partite hermitian two-faced pairs (\cite{10}, \cite{7}, \cite{5}, \cite{11}, \cite{8}).

The present paper has three more sections besides the introduction. Section one contains preliminaries. Section two derives the main technical result, the bi-free extremal convolution for the distributions of hermitian bi-partite pairs in the case the variables are projections. The third section gives the formula for the extremal bi-free convolution in the general case of hermitian variables and the definitions of bi-free max-stable and max-infinitely-divisible laws.

\section{Preliminaries}
\label{sec1}

\subsection{}
\label{sec1.1} Throughout the preliminaries $(M,\tau)$ will denote a $W^*$-\linebreak probability space, that is $M$ is a von~Neumann algebra and $\tau$ a normal state. If ${\mathcal A}$ is a $C^*$-algebra by Proj$(A)$ we shall denote the hermitian projections $P = P^* = P^2 \in {\mathcal A}$. If ${\mathcal A} = N$ a von~Neumann algebra and $P,Q \in \mbox{Proj}(N)$ by $P \vee Q$ and $P \wedge Q$ we denote the least projection $\ge$ both $P$ and $Q$ and, respectively, the largest projection $\le$ both $P$ and $Q$. If ${\mathcal H}$ is the Hilbert space on which $N$ acts then $P \vee Q$ is the orthogonal projection onto $\overline{P{\mathcal H} + Q{\mathcal H}}$, while $P \wedge Q$ is the orthogonal projection onto $P{\mathcal H} \cap Q{\mathcal H}$.

\subsection{}
\label{sec1.2} If $M_h = \{m \in M \mid m = m^*\}$, we recall that the spectral order (\cite{1}; see also \cite{2}) on $M_h$ is defined by $a \prec b$ if the spectral projections satisfy
\[
E(a;[t,\infty)) \le E(b;[t,\infty))
\]
for all $t \in {\mathbb R}$. Clearly this extends to self-adjoint operator affiliated with $M$, since $E(a;[t,\infty)) \le E(b;[t,\infty))$ makes sense also under this more general assumption. If $a,b \in M_h$ then $a \vee b$ and $a \wedge b$ are defined by
\[
E(a \vee b;(t,\infty)) = E(a;(t,\infty)) \vee E(b;(t,\infty))
\]
and
\[
E(a \wedge b;[t,\infty)) = E(a;[t,\infty)) \wedge E(b;[t,\infty)).
\]
These definitions clearly work also more generally for affiliated self-adjoint operators. In essence the operators $a,b$ are replaced by the right-continuous decreasing family of projections $R \ni t \to E(a,(t,\infty)) \in \mbox{Proj}(M)$ and, respectively, by the left-continuous decreasing family of projections ${\mathbb R} \ni t \to E(a;[t,\infty)) \in \mbox{Proj}(M)$.

\subsection{}
\label{sec1.3} A basic fact underlying free extreme values is the following:

\bigskip
\noindent
{\bf Lemma 1.1.} {\em If $P,Q \in \mbox{Proj}(M)$ are freely independent in $(M,\tau)$, then
\[
\tau(P \vee Q) = \min(\tau(P)+\tau(Q),1)
\]
and
\[
\tau(P \wedge Q) = \max(0,\tau(P)+\tau(Q)-1).
\]
}

\bigskip
This is a well-known fact. For references see \cite{2} where this is Lemma~2.1 and see \cite{12}, \cite{13} for computations. Remark that it is not necessary to require that $\tau$ be tracial, since its restrictions to the von~Neumann algebra generated by two free hermitian operators is always tracial.

We also recall from \cite{2} the definitions of extremal free convolutions for probability measures on ${\mathbb R}$. If $\mu$ is a probability measures on ${\mathbb R}$, its distribution function is $F(t) = \mu((-\infty,t])$. If $\mu$ and $\nu$ are probability measure on ${\mathbb R}$ with distribution functions $F(t)$ and $G(t)$, then $\mu \boxed{\vee} \nu$ and $\mu\boxed{\wedge}\nu$ are defined via their distribution functions $H(t) = \max(0,F(t)+G(t)-1)$ and, respectively, $K(t) = \min(F(t)+G(t),1)$. If $\mu$ and $\nu$ are the distributions of $a,b \in M_h$ with respect to $\tau$, then $\mu\boxed{\vee}\nu$ and $\mu\boxed{\wedge}\nu$ are the distributions of $a \vee b$ and $a \wedge b$. It is also convenient to have corresponding operations on distribution functions of probability measures on ${\mathbb R}$. If $F,G$ are two such distribution functions, then $F\boxed{\vee}G = (F+G-1)_+$ and $F\boxed{\wedge}G = \min(F+G,1)$.

\subsection{}
\label{sec1.4} We conclude the section of preliminaries by recalling some basics about the free $R$-transform and the partial bi-free $R$-transform.

If $a \in ({\mathcal A},\varphi)$ is a non-commutative random-variable in a Banach-algebra probability space, let $G_a(z) = \varphi((z1-a)^{-1})$ be the Green-function, or Cauchy-transform of the distribution of $a$, which is a holomorphic function in a neighborhood of $\infty \in {\mathbb C} \cup \{\infty\}$ on the Riemann sphere. Then $K_a(z)$ defined in a neighborhood of $0 \in {\mathbb C}$ and taking values in ${\mathbb C} \cup \{\infty\}$ is the inverse of $G_a$, that is $K_a(0) = \infty$ and $G_a(K_a(z)) = z$. Further $R_a(z) = K_a(z) - z^{-1}$ is defined in a neighborhood of $0$. If $a$ and $b$ are free, then $R_{a+b}(z) = R_a(z) + R_b(z)$ in a neighborhood of $0$.

If $(a,b)$ is a two-faced bi-free pair in $({\mathcal A},\varphi)$ we consider
\[
G_{a,b}(z,w) = \varphi((z1-a)^{-1}(w1-b)^{-1})
\]
defined in a neighborhood of $(\infty,\infty) \in ({\mathbb C} \cup \{\infty\})^2$. The reduced partial bi-free $R$-transform of $(a,b)$ is
\[
\widetilde{\mathcal R}_{a,b}(z,w) = 1 - \frac {zw}{G_{a,b}(K_a(z),K_b(w))}
\]
defined in a neighborhood of $(0,0) \in {\mathbb C}^2$. If $(a,b)$ and $(c,d)$ are bi-free in $({\mathcal A},\varphi)$, then
\[
\widetilde{\mathcal R}_{a+c,b+d}(z,w) = \widetilde{\mathcal R}_{a,b}(z,w) + \widetilde{\mathcal R}_{c,d}(z,w)
\]
in a neighborhood of $(0,0) \in {\mathbb C}^2$ (see \cite{10}).

If $(a,b)$ are commuting hermitian operators in a $C^*$-probability space $({\mathcal A},\varphi)$, then the joint distribution of $(a,b)$ is given by a probability measure on ${\mathbb R}^2$ with compact support $\mu_{a,b}$ and
\[
G_{a,b}(z,w) = \iint (z-x)^{-1}(w-y)^{-1}d\mu_{a,b}(x,y).
\]
The measures $\mu_{a,b}$, $\mu_{c,d}$ and $\mu_{a+c,b+d}$ when $(a,b)$ and $(c,d)$ are bi-free in $({\mathcal A},\varphi)$ are related by additive bi-free convolution
\[
\mu_{a,b} \boxplus\boxplus \mu_{c,d} = \mu_{a+c,b+d}.
\]

\section{Two-faced pairs of commuting projections}
\label{sec2}

In this section we shall compute the bi-free extremal convolution in the case of the distributions of commuting projections. This is the bi-free generalization of the free probability result in Lemma~1.1. We begin with a sequence of lemmas.

\bigskip
\noindent
{\bf Lemma 2.1.} {\em Let $({\mathcal A},\varphi)$ be a $C^*$-probability space and let $P = P^* = P^2 \in {\mathcal A}$ and $\varphi(P) = p$. Then we have
\[
G_P(z) = \frac {p}{z-1} + \frac {1-p}{z} = \frac {z+p-1}{z(z-1)}
\]
and
\[
zK_P(K_P-1) = K_P+p-1.
\]
}

\bigskip
\noindent
{\bf Lemma 2.2.} {\em Let $({\mathcal A},\varphi)$ be a $C^*$-probability space and let $(P,Q)$ be a two-faced pair in ${\mathcal A}$, so that $P = P^*$, $Q = Q^*$, $P = P^2$, $Q = Q^2$, $[P,Q] = 0$, $\varphi(P) = p$, $\varphi(Q) = q$, $\varphi(PQ) = r$. Then we have:
\[
G_{P,Q}(z,w) = \frac {(z+p-1)(w+q-1)+(r-pq)}{zw(z-1)(w-1)}\,.
\]
}

\bigskip
\noindent
{\bf Lemma 2.3.} {\em Under the same assumptions as in Lemma~$2.2$, we have:
\[
G_{P,Q}(K_P,K_Q) = \frac {zw((K_P+p-1)(K_Q+q-1)+(r-pq))}{(K_P+p-1)(K_Q+q-1)}\,.
\]
}

\bigskip
\noindent
{\bf Lemma 2.4.} {\em Under the same assumptions as in Lemma~$2.2$, we have:
\[
\widetilde{\mathcal R}_{P,Q}(z,w) = \frac {r-pq}{(K_P(z)+p-1)(K_Q(w)+q-1)+r-pq}\,.
\]
}

\bigskip
The proofs of the preceding four lemmata are straightforward computations and will be omitted.

\bigskip
\noindent
{\bf Lemma 2.5.} {\em Let $\mu$ be a probability measure on $[0,2]^2 \subset {\mathbb R}^2$ and let
\[
G(z,w) = \iint (z-x)^{-1}(w-y)^{-1}d\mu(x,y).
\]
Let further $z_n \in (2,\infty)$, $w_n \in (2,\infty)$ be such that $z_n \to 2$, $w_n \to 2$ as $n \to \infty$. Then we have:
\[
\lim_{n \to \infty}(z_n-2)(w_n-2)G(z_n,w_n) = \mu(\{(2,2)\}).
\]
}

\bigskip
\noindent
{\bf {\em Proof.}} Let $F_n(x,y) = (z_n-2)(z_n-x)^{-1}(w_n-2)(w_n-y)^{-1}$ where $(x,y) \in [0,2]^2$. Then $|F_n| \le 1$ for $(x,y) \in [0,2]^2$ and $F_n$ converges pointwise to the indicator function of $\{(2,2)\}$. By Lebesgue's dominated convergence theorem we have
\[
\lim_{n \to \infty} \iint F_n(x,y)d\mu(x,y) = \mu(\{(2,2)\})
\]
which is what we wanted to prove.\qed

\bigskip
\noindent
{\bf Lemma 2.6.} {\em In a $C^*$-probability space $({\mathcal A},\varphi)$ let $(P,Q)$ and $(P',Q')$ be bi-free two-faced pairs so that $P = P^* = P^2$, $Q = Q^* = Q^2$, $P' = P'{}^* = P'{}^2$, $Q' = Q'{}^* = Q'{}^2$, $[P,Q] = [P',Q'] = [P,Q'] = [P',Q] = 0$ and $\varphi(P) = p$,  $\varphi(Q) = q$, $\varphi(PQ) = r$, $\varphi(P') = p'$, $\varphi(Q') = q'$, $\varphi(P'Q') = r'$ and let $\delta = r-pq$, $\delta' = r'-p'q'$. Then for $(z,w)$ in some neighborhood of $(0,0) \in {\mathbb C}^2$ we have
\[
\begin{aligned}
&G_{P+P',Q+Q'}(K_{P+P'}(z),K_{Q+Q'}(w)) = \\
&= zw(1-(1+\delta^{-1}(K_P(z)+p-1)(K_Q(w)+q-1))^{-1} \\
&\quad - (1+\delta'{}^{-1}(K_{P'}(z)+p'-1)(K_{Q'}(w)+q'-1))^{-1})^{-1}.
\end{aligned}
\]
In case $\delta = 0$ we set here $(1+\delta^{-1}(K_P(z)+p-1)(K_Q(w)+q-1))^{-1} = 0$ and we adopt also a similar rule if $\delta' = 0$.
}

\bigskip
\noindent
{\bf {\em Proof.}} We have
\[
\begin{split}
G_{P+P',Q+Q'}(K_{P+P'}(z),K_{Q+Q'}(w)) &= zw(1-\widetilde{\mathcal R}_{P+P',Q+Q'}(z,w))^{-1} \\
&= zw(1-\widetilde{\mathcal R}_{P,Q}(z,w) \\
&\quad - \widetilde{\mathcal R}_{P',Q'}(z,w))^{-1}.
\end{split}
\]
Note that if $\delta = 0$, $P$ and $Q$ are classically independent, so that $\widetilde{\mathcal R}_{P,Q}(z,w) = 0$ and similarly $\widetilde{\mathcal R}_{P',Q'}(z,w) = 0$ if $\delta' = 0$. On the other hand, Lemma~2.4 gives that $\widetilde{\mathcal R}_{P,Q}(z,w) = (1 + \delta^{-1}(K_P(z)+p-1)(K_Q(w)+q-1))^{-1}$ if $\delta \ne 0$ and a similar fact for $\widetilde{\mathcal R}_{P',Q'}(z,w)$.
\qed

\bigskip
\noindent
{\bf Lemma 2.7.} {\em Under the same assumptions as in Lemma~$2.6$ we have for $(z,w)$ in some neighborhood of $(0,0) \in {\mathbb C}^2$ that
\[
\begin{split}
&(K_{P+P'}(z)-2)(K_{Q+Q'}(w)-2)G_{P+P',Q+Q'}(K_{P+P'}(z),K_{Q+Q'}(w)) \\
&= (K_P(z)+K_{P'}(z)-z^{-1}-2)(K_Q(w)+K_{Q'}(w)-w^{-1}-2)\\
&\qquad zw(1-(1+\delta^{-1}(K_P(z)+p-1)(K_Q(w)+q-1))^{-1} \\
&\quad - (1+\delta'{}^{-1}(K_{P'}(z)+p'-1)(K_{Q'}(w)+q'-1))^{-1})^{-1},
\end{split}
\]
this being an equality of holomorphic functions.
}

\bigskip
\noindent
{\bf {\em Proof.}} This follows from the preceding lemma after multiplication with
\[
\begin{split}
&(K_{P+P'}(z)-2)(K_{Q+Q'}(w)-2) = \\
&\quad(K_P(z)+K_{P'}(z)-z^{-1}-2)(K_Q(w)+K_{Q'}(w)-w^{-1}-2).
\end{split}
\]
Since the conclusion of Lemma~2.6 was actually an equality of germs of holomorphic functions near $(0,0) \in {\mathbb C}^2$ it might seem that there may be a problem with infinities when $z$ or $w$ is $0$. It is easily seen looking at the right-hand side that this is not the case since $(K_P(z)+K_{P'}(z)-z^{-1}-2)z$ is holomorphic in a neighborhood of $z=0$ and $(K_Q(w)+K_{Q'}(w)-w^{-1}-2)w$ is holomorphic in a neighborhood of $w=0$, while
\[
\begin{split}
&(1-(1+\delta^{-1}(K_P(z)+p-1)(K_Q(w)+q-1))^{-1} \\
&\quad - (1+\delta'{}^{-1}(K_{P'}(z)+p'-1)(K_{Q'}(w)+q'-1))^{-1})
\end{split}
\]
is holomorphic in a neighborhood of $(0,0) \in {\mathbb C}^2$ and $\ne 0$.
\qed

\bigskip
\noindent
{\bf Lemma 2.8.} {\em Assuming $p > 0$, $q > 0$, $p' > 0$, $q' > 0$, the equality which is the conclusion of Lemma~$2.7$ for $(z,w)$ in a neighborhood of $(0,0) \in {\mathbb C}^2$ extends analytically to $(z,w)$ in a neighborhood of $[0,\infty)^2 \subset {\mathbb C}^2$.
}

\bigskip
\noindent
{\bf {\em Proof.}} Let $t_0$ and $s_0$ be the least upper bounds of the supports of the probability measures $\mu_{P+P'}$ and $\mu_{Q+Q'}$ on ${\mathbb R}$. Then $G_{P+P'}(z)$ on $(t_0,\infty]$ and $G_{Q+Q'}(w)$ on $(s_0,\infty]$ are strictly decreasing taking the values $[0,\infty)$ so that $K_{P+P'}(z)$ and $K_{Q+Q'}(w)$ have analytic continuations along $[0,\infty)$ to a neighborhood of $[0,\infty)$ (the functions are viewed as taking values in the Riemann sphere ${\mathbb C} \cup \{\infty\}$). This implies the analytic extension of
\[
\begin{split}
&zw(K_P(z)+K_{P'}(z)-z^{-1}-2)(K_Q(w)+K_{Q'}(w)-w^{-1}-2) \\
&= zw(K_{P+P'}(z)-2)(K_{Q+Q'}(w)-2)
\end{split}
\]
as a holomorphic function, taking values in ${\mathbb C}$, to a neighborhood of $[0,\infty)^2 \subset {\mathbb C}^2$. On the other hand, similarly, since $p > 0$, $p' > 0$, $q > 0$, $q' > 0$ the functions $K_P(z)$, $K_{P'}(z)$, $K_Q(w)$, $K_{Q'}(w)$ have analytic continuations to a neighborhood of $[0,\infty)$ taking values for $z,w \in [0,\infty)$ in $(1,\infty]$. If $z,w \in [0,\infty)$, $(1+\delta^{-1}(K_P(z)+p-1)(K_Q(w)+q-1))^{-1}$ is well-defined when $\delta \ne 0$ since $(K_P(z)+p-1)(K_Q(w)+q-1) > pq > 0$ while either $\delta^{-1} > 0$ or $0 > \delta^{-1} \ge -p^{-1}q^{-1}$ so that the quantity to be inverted is $> 0$. Similar reasoning takes care of the term involving $\delta'$. From $[0,\infty)^2$ the extension goes over to a neighborhood.

For the analytic extension of
\[
G_{P+P',Q+Q'}(K_{P+P'}(z),K_{Q+Q'}(w))
\]
it suffices to remark that $G_{P+P',Q+Q'}(z,w)$ is analytic in a neighborhood of $(t_0,\infty] \times (s_0,\infty] \subset ({\mathbb C} \cup \{\infty\})^2$ and $K_{P+P'}(z)$, $K_{Q+Q'}(w)$ on $[0,\infty)$ take values in $(t_0,\infty]$ and $(s_0,\infty]$, respectively. That
\[
\begin{split}
&(1-(1+\delta^{-1}(K_P(z)+p-1)(K_Q(w)+q-1))^{-1} \\
&\quad - (1+\delta'{}^{-1}(K_{P'}(z)+p'-1)(K_{Q'}(w)+q'-1))^{-1})
\end{split}
\]
is $\ne 0$ for $(z,w) \in [0,\infty)^2$ follows from the fact the left-hand side is finite and
\[
\begin{split}
&zw(K_P(z)+K_{P'}(z)-z^{-1}-2)(K_Q(w)+K_{Q'}(w)-w^{-1}-2) \\
&= zw(K_{P+P'}(z)-2)(K_{Q+Q'}(w)-2) \ne 0.
\end{split}
\]
\qed

\bigskip
\noindent
{\bf Lemma 2.9.} {\em If $({\mathcal A},\varphi)$ is a $W^*$-probability space and $P = P^* = P^2 \in {\mathcal A}$, $P' = P'{}^* = P'{}^2 \in {\mathcal A}$ then $P \wedge P' = E(P+P',\{2\})$.
}

\bigskip
The preceding lemma is a well-known fact.

\bigskip
\noindent
{\bf Lemma 2.10.} {\em Let $({\mathcal A},\varphi)$ be a $W^*$-probability space. Then under the assumptions of Lemma~$2.6$, if $p + p' - 1 > 0$ and $q + q' - 1 > 0$ we have
\[
\begin{split}
&\varphi((P \wedge P')(Q \wedge Q')) \\
&= (p+p'-1)(q+q'-1)(1 - (1 + \delta^{-1}pq)^{-1} - (1 + \delta'{}^{-1}p'q')^{-1})^{-1}.
\end{split}
\]
Here in case $\delta = 0$ we set $(1+\delta^{-1}pq)^{-1} = 0$ and adopt a similar rule if $\delta' = 0$. In case $p+p'-1 \le 0$ or $q-q'-1 \le 0$ we have $\varphi((P \wedge P')(Q \wedge Q')) = 0$.
}

\bigskip
\noindent
{\bf {\em Proof.}} As recorded in Lemma~1.1, we have
\[
\varphi(P \wedge P') = (p+p'-1)_+,\ \varphi(Q\wedge Q') = (q+q'-1)_+
\]
after observing that $P$ and $P'$ being free, the restriction of $\varphi$ to the algebra generated by $P$ and $P'$ is a trace and a similar fact for $Q$ and $Q'$. Clearly, if $\varphi(P \wedge P') = 0$ or $\varphi(Q \wedge Q') = 0$, we must have also $\varphi((P \wedge P')(Q \wedge Q')) = 0$ since $[P \wedge P',Q \wedge Q'] = 0$ and $0 \le (P \wedge P')(Q \wedge Q') \le P \wedge P'$. Thus we are left with proving the lemma when $p+p'-1 > 0$ and $q+q'-1 > 0$.

Thus $\mu_{P+P'}(\{2\}) = p+p'-1 > 0$, $\mu_{Q+Q'}(\{2\}) = q+q'-1 > 0$ and $\mbox{supp } \mu_{P+P'} \subset [0,2]$, $\mbox{supp } \mu_{Q+Q'} \subset [0,2]$. By considerations along the lines in the proof of Lemma~2.8, if $t_n \in (0,\infty)$, $t_n \uparrow \infty$ then $K_P(t_n) \downarrow 1$, $K_{P'}(t_n) \downarrow 1$, $K_Q(t_n) \downarrow 1$, $K_{Q'}(t_n) \downarrow 1$, $K_{P+P'}(t_n) \downarrow 2$, $K_{Q+Q'}(t_n) \downarrow 2$. Taking $z = w = t_n$ in the equality in Lemma~2.7 extended according to Lemma~2.8, we get that the limit of the left-hand side in view of Lemma~2.5 is
\[
\begin{split}
\mu_{P+P',Q+Q'}(\{(2,2)\}) &= \varphi(E(P+P',\{2\})E(Q+Q',\{2\})) \\
&= \varphi((P \wedge P')(Q \wedge Q')).
\end{split}
\]

On the other hand,
\[
\begin{split}
&(K_P(t_n)+K_{P'}(t_n)-t_n^{-1}-2)t_n \\
&= (K_P(t_n)-1)t_n + (K_{P'}(t_n)-1)t_n-1 \\
&= G_P(K_P(t_n))(K_P(t_n)-1)+G_{P'}(K_{P'}(t_n))(K_{P'}(t_n)-1)-1
\end{split}
\]
and this converges as $n \to \infty$, by the simpler analogue of Lemma~2.5 for Cauchy transforms in one variable, to $p+p'-1$. Similarly $(K_Q(t_n)+K_{Q'}(t_n)-t_n^{-1}-2)t_n$ converges to $q+q'-1$. On the other hand
\[
\begin{split}
&(1-(1+\delta^{-1}(K_P(t_n)+p-1)(K_Q(t_n)+q-1))^{-1} \\
&\quad - (1+\delta'{}^{-1}(K_{P'}(t_n)+p'-1)(K_{Q'}(t_n)+q'-1))^{-1})
\end{split}
\]
converges to $(1-(1+\delta^{-1}pq)^{-1}-(1+\delta'{}^{-1}p'q'))$ no matter whether $\delta$ and $\delta'$ are $\ne 0$ or $= 0$.
\qed

\bigskip
The last lemma after some simple algebraic work on the formulae will give the final result of the computations in this section, which we record as a theorem.

\bigskip
\noindent
{\bf Theorem 2.1.} {\em Let $({\mathcal A},\varphi)$ be a $W^*$-probability space and let $P = P^* = P^2 \in {\mathcal A}$, $Q = Q^* = Q^2 \in {\mathcal A}$, $P' = P'{}^* = P'{}^2 \in {\mathcal A}$, $Q' = Q'{}^* = Q'{}^2 \in {\mathcal A}$ be such that $[P,Q] = 0$, $[P',Q] = [P,Q'] = 0$, $[P',Q'] = 0$ and $(P,Q)$ and $(P',Q')$ are bi-free in $({\mathcal A},\varphi)$. Then we have $\varphi(P \wedge P') = (\varphi(P)+\varphi(P')-1)_+$, $\varphi(Q \wedge Q') = (\varphi(Q)+\varphi(Q')-1)_+$ and if $\varphi(P\wedge P') > 0$, $\varphi(Q \wedge Q') > 0$, $\varphi(PQ) > 0$, $\varphi(P'Q') > 0$ then we have
\[
\frac {\varphi(P \wedge P')\varphi(Q \wedge Q')}{\varphi((P \wedge P')(Q \wedge Q'))} = \frac {\varphi(P)\varphi(Q)}{\varphi(PQ)} + \frac {\varphi(P')\varphi(Q')}{\varphi(P'Q')} - 1.
\]
If any of the numbers $\varphi(P\wedge P')$, $\varphi(Q \wedge Q')$, $\varphi(PQ)$, $\varphi(P'Q')$ is $0$, then $\varphi((P \wedge P')(Q \wedge Q')) = 0$.
}

\bigskip
\noindent
{\bf {\em Proof.}} The formulae for $\varphi(P \wedge P')$, $\varphi(Q \wedge Q')$ are not new (see Lemma~1.1) and it is obvious that if any of $\varphi(P \wedge P')$, $\varphi(Q \wedge Q')$, $\varphi(PQ)$ $\varphi(P'Q')$ is $0$, then so is $\varphi((P \wedge P')(Q \wedge Q'))$. Thus using the notation $\varphi(P) = p$, $\varphi(P') = p'$, $\varphi(Q) = q$, $\varphi(Q') = q'$, $\varphi(PQ) = r$, $\varphi(P'Q') = r'$ which we used in the lemmata, we may assume $p + p' > 1$, $q + q' > 1$, $r > 0$, $r' > 0$. Turning to the result of Lemma~2.10, remark that if $\delta = r - pq \ne 0$ then
\[
(1+\delta^{-1}pq)^{-1} = \delta(\delta+pq)^{-1} = (r-pq)r^{-1} = 1 - pqr^{-1}.
\]
If $\delta = 0$, then $r = pq$ and $1 - pqr^{-1} = 0$ which is in agreement with the rule that $(1+\delta^{-1}pq)^{-1}$ is $0$ if $\delta = 0$. A similar remark about $\delta'$. Hence the right-hand side of the formula in Lemma~2.10 is
\[
\varphi(P \wedge P')\varphi(Q \wedge Q')(pqr^{-1} + p'q'r'{}^{-1} - 1)^{-1}
\]
so that the formula gives
\[
\begin{split}
\frac {\varphi(P \wedge P')\varphi(Q \wedge Q')}{\varphi((P \wedge P')(Q \wedge Q'))} 
&= pqr^{-1} + p'q'r'{}^{-1} -1 \\
&= \frac {\varphi(P)\varphi(Q)}{\varphi(PQ)} + \frac {\varphi(P')\varphi(Q')}{\varphi(P'Q')} - 1.
\end{split}
\]
\qed

\section{Bi-free max-convolution in the plane}
\label{sec3}

In this section $({\mathcal A},\varphi)$ will be a von~Neumann algebra with a normal state $\varphi$. If $(a,b)$ is a bi-partite hermitian two-faced pair in $({\mathcal A},\varphi)$, that is a pair of commuting hermitian operators $a,b \in {\mathcal A}$, let $E(a,b;\omega)$ denote its spectral measure where $\omega \subset {\mathbb R}^2$ is a Borel set and let $\mu_{a,b}(\omega) = \varphi(E(a,b;\omega))$ be the probability measure on ${\mathbb R}^2$ which is the distribution of $(a,b)$. Let further $F_{a,b}(s,t) = \mu_{a,b}((-\infty,s] \times (-\infty,t])$ be the distribution function of the measure $\mu_{a,b}$. We recall that such functions $F(s,t)$ are such that $s_1 \le s_2$, $t_1 \le t_2 \Rightarrow F(s_1,t_1) \le F(s_2,t_2)$, $s_n \downarrow s_0$, $t_n \downarrow t_0 \Rightarrow F(s_n,t_n) \downarrow F(s_0,t_0)$ as $n \to \infty$ and $s_1 \le s_2$, $t_1 \le t_2 \Rightarrow F(s_2,t_2) - F(s_1,t_2) - F(s_2,t_1) + F(s_1,t_1) \ge 0$. Moreover, since this is the distribution function of a probability measure with compact support, we have $0 \le F(s,t) \le 1$ and there is $C > 0$ so that $F(s,t) = 0$ if $\min(s,t) \le -C$ and $F(s,t) = 1$ if $\min(s,t) \ge C$. If we want to deal with probability measures without a condition of compact support, we will require that $F$ be defined on $[-\infty,\infty)^2$ and satisfy $F(s,t) = 0$ if $\min(s,t) = -\infty$ and $\underset{n \uparrow +\infty}{\lim}\ F(n,n) = 1$.

If $(a,b)$ and $(a',b')$ are bi-free bi-partite hermitian pairs, it's always possible to find a realization in a von~Neumann algebra $({\mathcal A},\varphi)$ of the joint distribution so that $[a,b'] = [a',b] = 0$. Note further that the joint distribution of $(a \vee a',b \vee b')$ does not depend on the realization, but only on the distributions $\mu_{a,b}$, $\mu_{a',b'}$ since
\[
\begin{split}
&E(a \vee a',b \vee b';(-\infty,s] \times (-\infty,t]) \\
&= E(a \vee a';(-\infty,s])E(b \vee b';(-\infty,t]) \\
&= (E(a;(-\infty,s]) \wedge E(a';(-\infty,s]))(E(b;(-\infty,t]) \wedge E(b';(-\infty,t]))
\end{split}
\]
and $\varphi(E(a \vee a',b \vee b';(-\infty,s] \times (-\infty,t]))$ can be computed using the results in Section~2 from the distributions of the bi-free two-faced pairs $(E(a;(-\infty,s],E(b;(-\infty,t]))$ and $(E(a';(-\infty,s]),E(b';(-\infty,t]))$.

\bigskip
\noindent
{\bf Definition 3.1.} If $F$ and $G$ are distribution functions of probability measures with compact support on ${\mathbb R}^2$ we define their bi-free max-convolution or alternatively also called bi-free upper extremal convolution $H = F\boxed{\vee}\boxed{\vee}G$ to be such that if $F_j,G_j,H_j$ $(j = 1,2)$ are their distribution function marginals we have $H_j = F_j\boxed{\vee}G_j$ $(j = 1,2)$ and
\[
\frac {H_1(s)H_2(t)}{H(s,t)} = \frac {F_1(s)F_2(t)}{F(s,t)} + \frac {G_1(s)G_2(t)}{G(s,t)} - 1
\]
if $F(s,t) > 0$, $G(s,t) > 0$, $H_1(s) > 0$, $H_2(t) > 0$ and $H(s,t) = 0$ otherwise.

That the above gives a well-defined distribution function of a probability measure with compact support is a consequence of the discussion preceding the definition and of Theorem~2.1. Note also that to see that the distribution function of $a \wedge a'$ is given by $F_a\boxed{\vee}F_{a'}$ it is not necessary to assume $a$ and $a'$ are in a tracial $W^*$-probability space, since the restriction of $\varphi$ to the weak closure of the $*$-algebra generated by $\{I,a,a'\}$ will be a tracial normal state. In essence, $\boxed{\vee}\boxed{\vee}$ gives the distribution of $(a \vee a',b \vee b')$ in the realizations of the joint distributions of $(a,a')$ and $(b,b')$ where the commutations $[a,b'] = [a',b] = 0$ hold.

The further remark is that actually $a,b,a',b'$ only appear here via their spectral scales $E(a;(-\infty,t])$ etc. and thus {\em the operations extend to affiliated unbounded self-adjoint operators and distributions of any probability measures on} ${\mathbb R}^2$.

Having defined $\boxed{\vee}\boxed{\vee}$ we can now define bi-free max-stable and bi-free max-infinitely divisible laws on ${\mathbb R}^2$.

\bigskip
\noindent
{\bf Definition 3.2.} A distribution function $F$ of a probability measure on ${\mathbb R}^2$ is bi-freely max-stable if there are $a_n,b_n,c_n,d_n \in {\mathbb R}$, $a_n > 0$, $c_n > 0$ so that
\[
(\underset{\mbox{$n$-fold}}{\underbrace{F\boxed{\vee}\boxed{\vee}F\boxed{\vee}\boxed{\vee}F\boxed{\vee}\boxed{\vee}\dots F}})(a_nx+b_n,c_ny+d_n) \to F(x,y)
\] 
as $n \to \infty$.

\bigskip
\noindent
{\bf Definition 3.3.} A distribution function $F$ of a probability measure on ${\mathbb R}^2$ is bi-freely max-infinitely-divisible if for each $n \in {\mathbb N}$ there is a distribution function $F_n$ so that
\[
\underset{\mbox{$n$-fold}}{\underbrace{F_n\boxed{\vee}\boxed{\vee}F_n\boxed{\vee}\boxed{\vee}\dots F_n}} = F.
\]

\bigskip
\noindent
{\bf Remark 3.1.} The definitions in this section show that, in the simplest bi-free case of bi-partite hermitian two-faced pairs that have distributions given by probability measures on ${\mathbb R}^2$, the basic extreme value questions about bi-free max-stable and bi-free max-infinitely are transformed by the operation $\boxed{\vee}\boxed{\vee}$ into ``classical'' questions. Clearly these questions are more difficult than univariate free extreme value questions (\cite{2}, \cite{3}). It is a natural question whether like in \cite{2}, where free max-stable laws were related to classical ``peaks over thresholds'', the ``classical'' questions to which bi-free extremes laws lead in this simplest case are also related to some classical extremes questions (\cite{4}, \cite{6}).

\end{document}